\documentclass[letterpaper, 10 pt, conference]{ieeeconf}  

\IEEEoverridecommandlockouts                              
\usepackage{graphicx}      
\usepackage{amsfonts}
\usepackage{amsmath}
\usepackage{amssymb}
\usepackage{makeidx}
\usepackage{epsfig}

\def\R{\mathcal{R}}
\def\P{\mathcal{P}}
\def\O{\mathcal{O}}

\def\M{\mathcal{M}}

\newtheorem{assumption}{Assumption}
\newtheorem{theorem}{Theorem}
\newtheorem{corollary}[theorem]{Corollary}
\newtheorem{definition}{Definition}
\newtheorem{example}{Example}
\newtheorem{lemma}[theorem]{Lemma}
\newtheorem{proposition}[theorem]{Proposition}

\title{\LARGE \bf Fault detection and isolation of malicious nodes\\in MIMO Multi-hop Control Networks}

\author{A. D'Innocenzo, M.D. Di Benedetto and F. Smarra
\thanks{The authors are with the Center of Excellence DEWS and with the Department of Information Engineering, Computer Science and Mathematics of the University of L'Aquila, Italy. The research leading to these results has received funding from the European Union Seventh Framework Programme [FP7/2007-2013] under grant agreement n°257462 HYCON2 Network of excellence}}

\begin{document}

\maketitle
\thispagestyle{empty}
\pagestyle{empty}

\begin{abstract}
A MIMO Multi-hop Control Network (MCN) consists of a MIMO LTI system where the communication between sensors, actuators and computational units is supported by a (wireless) multi-hop communication network, and data flow is performed using scheduling and routing of sensing and actuation data.
We provide necessary and sufficient conditions on the plant dynamics and on the communication protocol configuration such that the Fault Detection and Isolation (FDI) problem of failures and malicious attacks to communication nodes can be solved.
\end{abstract}

\section{Introduction} \label{secIntro}

\begin{figure*}[t]
\begin{center}
\includegraphics[width=1\textwidth]{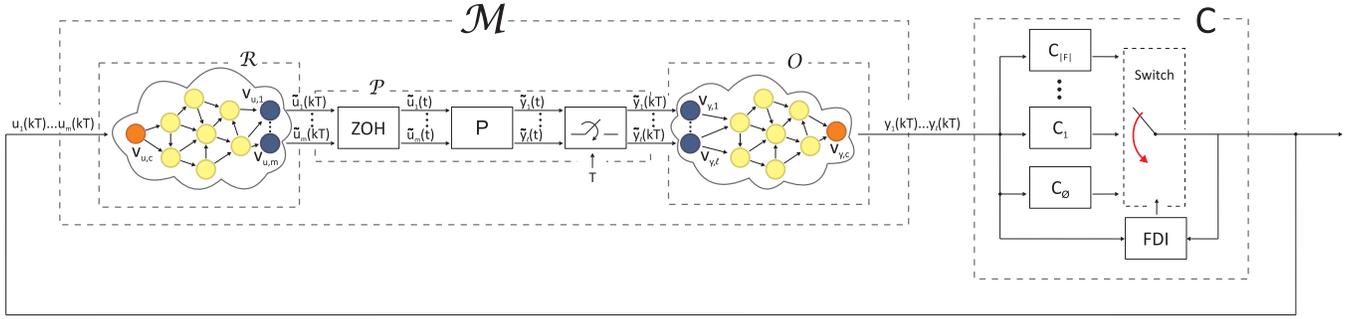}
\vspace{-0.4cm}
\caption{Control scheme of a Multi-hop Control Network subject to link failures.} \label{figHybridControlScheme}
\vspace{-0.8cm}
\end{center}
\end{figure*}

Wireless networked control systems are spatially distributed control systems where the communication between sensors, actuators, and computational units is supported by a wireless multi-hop communication network. The use of wireless Multi-hop Control Networks (MCNs) in industrial automation results in flexible architectures and generally reduces installation, debugging, diagnostic and maintenance costs with respect to wired networks. Although MCNs offer many advantages, co-design of the network configuration and of the control algorithm for a MCN requires addressing the joint dynamics of the plant and of the communication protocol.

Recently, a huge effort has been made in scientific research on Networked Control Systems (NCSs), see e.g. \cite{Astrom97j1},~\cite{Zhang2001},~\cite{Arzen06},~\cite{TabbaraTAC2007},~\cite{Hespanha2007},~\cite{MurrayTAC2009},~\cite{HeemelsTAC11} and references therein for a general overview. In general, the literature on NCSs addresses non--idealities (such as quantization errors, packets dropouts, variable sampling and delay and communication constraints) as aggregated network performance variables, losing irreversibly the dynamics introduced by scheduling and routing communication protocols. What is needed for modeling and analyzing control protocols on MCNs is an integrated framework for analysing/co-designing network topology, scheduling, routing and control. In~\cite{Andersson:CDC05}, a simulative environment of computer nodes and communication networks interacting with the continuous-time dynamics of the real world is presented. To the best of our knowledge, the first formal model of a Multi-hop Control Network has been presented in \cite{AlurRTAS09,AlurTAC11}, where a mathematical framework has been proposed that allows modeling the MAC layer (communication scheduling) and the Network layer (routing) of recently developed wireless industrial control protocols, such as WirelessHART and ISA-100.
In \cite{DInnocenzoTAC13} we extended the formalism proposed in \cite{AlurTAC11} by defining a MCN $\M$, that consists of a continuous-time SISO LTI plant $\P$ interconnected to a controller $C$ via two (wireless) multi-hop communication networks $G_\R$ (the controllability network) and $G_\O$ (the observability network), as illustrated in Figure \ref{figHybridControlScheme}, and by modeling redundancy in data communication - i.e. sending actuation/sensing data via multiple paths and then merging these components according to a weight function. This approach, which can be interpreted as a form of network coding at the level of the application layer of the ISO/OSI protocol stack, is called \emph{multi-path routing} (or \emph{flooding}, in the communication scientific community) and aims at enabling the detection and isolation of node failures and malicious intrusions, which cannot be done using single-path routing or strategies that use timestamps to discard redundant packets. It is well known that redundancy can also render the system fault-tolerant with respect to node failures and mitigate the effect of packet losses. We remark that, as illustrated in \cite{NeanderISIE2011}, the implementation of multi-path routing in a Wireless HART device only requires a minor change which retains backward compatibility with standard devices.

\textbf{Paper contribution:} Because of wireless networking, a MCN can be subject to failures and/or malicious attacks. In \cite{DInnocenzoTAC13} we addressed and solved the problem of designing a set of controllers and the communication protocol parameters of a SISO MCN, so that it is possible to detect and isolate the faulty nodes of the controllability and observability networks and apply an appropriate controller to stabilize the system, as depicted in Figure \ref{figHybridControlScheme}. In \cite{SmarraNecSys2012} and in this paper we extend such investigation to MIMO MCNs: in particular, while in \cite{SmarraNecSys2012} we extended the MCN formalism to MIMO LTI plants and developed a procedure to guarantee the existence of a stabilizing controller $C_{i}$ for any node failure, in this paper we provide conditions on network topology, scheduling and routing that enable detection and isolation of node failures.

The extension with respect to the results in \cite{DInnocenzoTAC13} is not trivial: indeed in the MIMO case the geometric approach exploited in \cite{DInnocenzoTAC13} does not easily provide a relation between the conditions that enable FDI of node failures and the network topology, scheduling and routing. To overcome this issue we exploit formalism and FDI methodologies of \emph{structured systems} \cite{DionAutomatica2003}. This methodology leverages on the classical observer-based results in \cite{MassoumniaTAC89}, where the system with failures is required to be left-invertible: therefore our results imply that \emph{almost any} failure signal can be detected and isolated. In order to detect and isolate \emph{any} non-zero failure signal we should require input-observability. Indeed a carefully chosen (malicious) failure may hide in the zero dynamics, being therefore undetectable: for example in {\cite{SundaramTAC11,PasqualettiTAC12}} the role of invariant zeros in detecting malicious attacks has

\textbf{Related work:} As can be inferred from the recent survey \cite{Gupta10}, fault tolerant control and fault diagnosis is one of the main issues addressed in the research on NCSs. However, most of the existing literature does not consider the effect of the communication protocol introduced by a Multi-hop Control Network. In \cite{CommaultTAC2007}, a procedure to minimize the number and cost of additional sensors, required to solve the FDI problem for \emph{structured systems}, is presented. In \cite{PappasCDC2010a} the design of an intrusion detection system is presented for a Wireless Control Network: our work differs from such results in the following four aspects. (1) In \cite{PappasCDC2011} the wireless network is an autonomous system where the network \emph{itself} acts as a decentralized controller, while in our model the wireless network \emph{transfers} sensing and actuation data between a plant and a centralized controller, namely it acts as a \emph{relay} network. This modeling choice is motivated by the fact that WirelessHART and ISA-100 are designed for control loops where a centralized controller exploits a wireless network to \emph{relay} sensing and actuation data, which is often a forced choice in industrial environments. (2) As a consequence of the previous issue, in \cite{PappasCDC2011} FDI is performed only exploiting the output signal from a subset of communication nodes, while we can exploit the input and output signals of the centralized controller. (3) In our model we take into account the effect of the scheduling ordering of the node transmissions in the sensing and actuation data relay, which provides a more accurate modeling of the effect of scheduling on the closed loop dynamics. Indeed, the conditions we derive in Section~\ref{secFDIMCN} show that, in order to guarantee FDI of faulty nodes, the link scheduling order is irrelevant: this is an interesting result because, as widely discussed in \cite{D'InnocenzoCASE2009,DInnocenzoTAC13,SmarraCDC12}, it strongly reduces the scheduling period, avoids the necessity of on-the-fly scheduling re-definition, and always guarantees the existence of an admissible scheduling when multiple loops exploit the same communication network. (4) In \cite{PappasCDC2010a} FDI is performed by on-the-fly testing the rank of a number of matrices which is a combinatorial function of the number of communication links, while our method only requires to apply a logic operator to a number of Luenberger observers which is at most equal to the number of communication nodes. To the best of our knowledge, our work is pioneering in addressing FDI for a MCN that implements standardized communication protocols. An extended version of this paper can be found on ArXiv.

\textbf{Notation:} We will denote by $\mathbb N$ and $\mathbb R$ respectively the sets of natural and real numbers.
Given $n \in \mathbb N$, we denote by $\bold{n}$ the set $\bold{n} \doteq \{1,2,\ldots,n\}$.
We denote by $\textbf{0}_{n \times m}$ the matrix of zeros with $n$ rows and $m$ columns and by $\textbf{I}_{n}$ the identity matrix of dimension $n$.
Given a finite set $A$ and a subset $B \subseteq A$, we define $|A|$ and $|B|$ their cardinality, $A \setminus B$ the difference set and $2^A$ the power set.
We denote by $diag(F_1(z), \ldots, F_n(z))$ the $n \times n$ diagonal transfer function matrix whose diagonal consists of the scalar transfer functions $F_1(z),\ldots,F_n(z)$.
Given a directed graph $(\mathcal V, \mathcal E)$ we define \emph{path} an alternating sequence of vertices and edges. A path is defined \emph{simple} if no vertices are repeated. We define a set of paths \emph{vertex disjoint} if each two of them consist of disjoint sets of vertices. We call a set of $r$ disjoint and simple paths from a set $\mathcal{V}_1 \subseteq \mathcal{V}$ to a set $\mathcal{V}_2 \subseteq \mathcal{V}$ an $r$\emph{-linking} from $\mathcal{V}_1$ to $\mathcal{V}_2$. We denote by $\biguplus$ the disjoint union operator among directed graphs. For a formal definition of further graph properties and notations (e.g. weakly connected graph, weakly connected component, bridge nodes etc.) the reader is referred to \cite{West01}.


\section{MCN model} \label{secModelingNominalMCNs}

\begin{figure}[ht]
\begin{center}
\vspace{-0.1cm}
\includegraphics[width=0.3\textwidth]{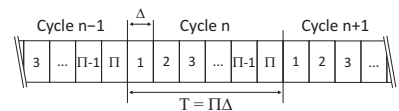}
\vspace{-0.3cm}
\caption{Time-slotted structure of frames.}\label{frame}
\vspace{-0.4cm}
\end{center}
\end{figure}

We propose a mathematical framework for modeling wireless multi-hop communication networks that implement time-triggered protocols such as WirelessHART and ISA-100. In these standards the access to the shared communication channel is specified as follows: time is divided into slots of fixed duration $\Delta$ and groups of $\Pi$ time slots are called frames of duration $T = \Pi \Delta$ (see Figure \ref{frame}). For each frame, a communication scheduling allows each node to transmit data only in a specified time slot. The scheduling is periodic with period $\Pi$, i.e. it is repeated in all frames.
\begin{definition}\label{defMCN}
A MIMO Multi-hop Control Network is a tuple $\M = ( \P, G, W, \eta, \Delta)$ where:

$\P$ is a continuous-time MIMO LTI system, with $n$, $m$ and $\ell$ respectively the dimensions of the internal state, input and output spaces.

$G = (G_{\R}, G_{\O})$. $G_{\R} = (V_{\R},E_{\R})$ is a directed graph, where the vertices correspond to the communication nodes of the network and an edge from $v$ to $v'$ means that node $v'$ can receive messages transmitted by node $v$ through the wireless communication link $(v,v')$. We denote by $v_{u,c}$ the special node of $V_{\R}$ that corresponds to the controller and by $v_{u,i} \in V_{\R}$, $i \in \bold{m}$, the special nodes that correspond to the actuators of the input components. $G_{\O} = (V_{\O},E_{\O})$ is defined similarly to $G_{\R}$. We denote by $v_{y,c}$ the special node of $V_{\O}$ that corresponds to the controller and by $v_{y,i} \in V_{\O}$, $i \in \boldsymbol{\ell}$, the special nodes that correspond to the sensors of the outputs $y_i$, $i \in \boldsymbol{\ell}$.

$W = (W_{\R}, W_{\O})$. $W_{\R} = \{W_{\R_i}\}_{i \in \bold{m}}$, where $W_{\R_i} : E_{\R} \to \mathbb R$ is a weight function for the $i$-th input component that associates to each link a real constant. $W_{\O} = \{W_{\O_i}\}_{i \in \boldsymbol{\ell}}$ is defined similarly to $W_{\R}$.

$\eta = (\eta_{\R}, \eta_{\O})$. $\eta_{\R} = \{\eta_{\R_i}\}_{i \in \bold{m}}$, where $\eta_{\R_i} \colon \mathbb N \to 2^{E_{\R}}$ is the controllability scheduling function for the $i$-th input component that associates to each time slot of each frame a set of edges of the controllability radio connectivity graph $G_{\R}$. Since in this paper we only consider a periodic scheduling, that is repeated in all frames, we define the controllability scheduling functions by $\eta_{\R_i} \colon \{1, \ldots, \Pi\} \to 2^{E_{\R}}$. The integer constant $\Pi$ is the period of the controllability scheduling. The semantics of $\eta_{\R_i}$ is that $( v,v' ) \in \eta_{\R_i}(h)$ if at time slot $h$ of each frame the data associated to the $i$-th input component and contained in node $v$ is transmitted to the node $v'$, multiplied by the weight $W_{\R_i}(v,v')$. For any $\eta_{\R_i}$, we assume that each link can be scheduled only one time for each frame\footnote{This does not lead to loss of generality, since it is always possible to obtain an equivalent model that satisfies this constraint by appropriately splitting the nodes of the graph, as already illustrated in the memory slot graph definition of \cite{AlurTAC11}.}. $\eta_{\O} = \{\eta_{\O_i}\}_{i \in \boldsymbol{\ell}}$ is defined similarly to $\eta_{\R}.$\footnote{We remark that the scheduling period of $\eta_{\O}$ is the same of $\eta_{\R}$.}

$\Delta$ is the time slot duration. As a consequence, ${T} = \Pi \Delta$ is the frame duration.
\end{definition}
The main difference with respect to the MCN definition given in \cite{DiBenedettoIFAC11Stab} for SISO systems is that here the network needs to relay $m$ actuation signals and $\ell$ sensor signals: we assume that nodes routes data of each of these signals separately.

Designing a scheduling function induces a communication scheduling (namely the time slot when each node is allowed to transmit) and a multi-path routing (namely the set of paths that convey data from the input to the output of the connectivity graph) of the communication protocol. Since the scheduling function is periodic, the induced communication scheduling is periodic and the induced multi-path routing is static.

\begin{definition}
Given $G_{\R}$ and $\eta_{\R_i}$ we define $G_{\R}\left(\eta_{\R_i}(h)\right)$ the sub-graph of $G_{\R}$ induced by keeping the edges scheduled in the time slot $h$. We define $G_{\R}(\eta_{\R_i}) = \bigcup_{h = 1}^{\Pi} G_{\R}(\eta_{\R_i}(h))$ the sub-graph of $G_{\R}$ induced by keeping the union of edges scheduled during the whole frame.
\end{definition}
$G_{\R}(\eta_{\R_i})$ consists of a set of simple (routing) paths starting from $v_{u,c}$ and terminating in $v_{u,i}$. As a consequence $G_{\R}(\eta_{\R_i})$ is a directed, weakly connected and acyclic graph. The above definition can be given similarly for $G_{\O}$ and $\eta_{\O_i}$.


\begin{figure}[ht]
\begin{center}
\includegraphics[width=0.48\textwidth]{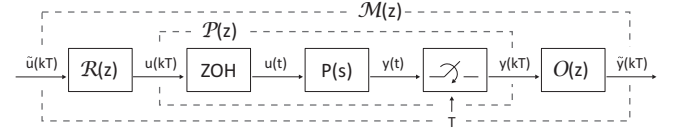}
\vspace{-0.7cm}
\caption{MCN interconnected system.} \label{MCNblocks}
\vspace{-0.5cm}
\end{center}
\end{figure}
\textbf{Nominal MCN:} The dynamics of a MCN $\M$ can be modeled by the interconnection of blocks as in Figure \ref{MCNblocks}. The block $\P$ is characterized by the discrete-time transfer function matrix $\P(z)$ obtained by discretizing the system $\P$ with sampling time ${T} = \Pi \Delta$. The block $\R$ models the dynamics introduced by the flow of the actuation data of all components of the control input $u$ via the communication network represented by $G_{\R}$. In order to define the dynamics of $\R$, we need to define the semantics of the data flow through the network induced by the scheduling and the weight functions. We assume that each communication node, when scheduled to transmit by $\eta_{\R_i}$, computes a linear combination of the data associated to the input component $u_i$ and received from the incoming links according to the weight function $W_{\R_i}$. This linear combination is then transmitted via the outgoing scheduled links. As in \cite{SmarraNecSys2012} the input/output behavior of $u_i(kT)$ with respect to $\tilde u_i(kT)$ can be formalized, for any $i \in \bold{m}$, by the following transfer function:
\begin{equation}\label{eqPropGz}
\R_i(z) \doteq \frac{\widetilde{U}_i(z)}{U_i(z)} = \sum\limits_{d = 1}^{D_{\R_i}} \frac{\gamma_{\R_i}(d)}{z^d},
\end{equation}
where $D_{\R_i} \in \mathbb N$ is the maximum delay introduced by the (routing) paths of $G_{\R}(\eta_{\R_i})$ and $\forall d \in \boldsymbol{D_{\R_i}}$, $\gamma_{\R_i}(d) \in \mathbb R$, $\gamma_{\R_i}(D_{\R_i}) \neq 0$. The reader is referred to \cite{DInnocenzoTAC13,SmarraCDC12} for the formal definition of the coefficients $\gamma_{\R_i}(d), d \in \boldsymbol{D_{\R_i}}$, which depend on the weight function $W_\R$.

The block $\R$ is characterized by the transfer function matrix $\R(z)=diag(\R_1(z),\ldots,\R_m(z))$.
The same holds for the block $\O$. The dynamics of a MIMO MCN $\M$ can be modeled by the cascade of the transfer function matrices $\R(z)$, $\P(z)$ and $\O(z)$, thus $\M(z) = \O(z) \P(z) \R(z)$.


\textbf{Faulty MCN:} We assume that a failure or a malicious attack associated to a communication node $v \in V_\R$ can be modeled by a set of arbitrary signals $f_{v,i}(k)$, for any $i$ such that there exist $v' \in V_\R$, $h \in \boldsymbol{\Pi}$ with $(v,v') \in \eta_{\R_i}(h)$, each summed to the $i$-th input component routed via node $v$.
This general framework, as illustrated in \cite{DInnocenzoTAC13}, models several classes of failures (e.g. a node stops sending data or sends fake data) and malicious attacks (e.g. an arbitrary signal is injected, which overrides/sums to the original data). Following the same reasoning as in the definition of $\R_i(z)$, we can define the transfer function from $f_{v,i}(k)$ to $u_i(k)$ as follows:
\begin{equation}\label{eqPropFz}
\mathcal{F}_{v,i}(z) \doteq \frac{U_i(z)}{F_{v,i}(z)} = \sum\limits_{d = 1}^{D_{v,i}} \frac{\gamma_{v,i}(d)}{z^d},
\end{equation}
where $D_{v,i} \in \mathbb N$ is the maximum delay introduced by the (routing) paths from $v$ to the actuator node $v_{u,i}$ and
$\forall d \in \boldsymbol{D_{v,i}}$, $\gamma_{v,i}(d) \in \mathbb R$, with $\gamma_{v,i}(D_{v,i}) \neq 0$. By the properties of $G_{\R}(\eta_{\R_i})$ it follows that $\forall v \in V_{\R}, D_{v,i} \leq D_{\R_i}$. The same holds for the block $\O$.


\section{MCN structured model} \label{secModelingFaultyMCNs}

To state conditions for FDI of failures and malicious attacks on a node of a MCN we will exploit the formalism of \emph{structured systems} \cite{DionAutomatica2003}. Given a system characterized by a state space representation $S = (A, B, C, D)$ we can define the associated structured system by defining the matrices $S_\lambda = (A_\lambda, B_\lambda, C_\lambda, D_\lambda)$ so that each entry is either zero (if the corresponding entry in the original matrix is zero) or a free parameter (if the corresponding entry of the original matrix is non-zero). For instance, consider a system $S$ given by $A = [1, 2; 0, 0]$, $B = [0, 1]^{\top}$, $C = [1, 0]$: the corresponding structured system $S_\lambda$ is given by $A = [\lambda_1, \lambda_2; 0, 0]$, $B = [0, \lambda_3]^{\top}$, $C = [\lambda_4, 0]$, where the $\lambda_i$'s are free parameters. A structured system can also be represented by a directed graph $(V_{S_\lambda}, E_{S_\lambda})$ whose vertices correspond to the input, state and output variables, and with an edge between two vertices if there is a non-zero free parameter $\lambda_i$ relating the corresponding variables in the equations. The graph representation of the example above is given by $V_{S_\lambda} = \{u, x_1, x_2, y\}$, and $E_{S_\lambda} = \{(x_1,x_1),(x_2,x_1),(u,x_2),(x_1,y)\}$.

The following proposition formalizes the graph structured representation of the block $\R$ when a set of failure signals is applied to communication nodes.
\begin{proposition}
Given $G_{\R}$, $\eta_{\R}$ and a set of faulty nodes $\bar V \subseteq V_\R$ we define the graph structured representation $(V_{\R_\lambda}, E_{\R_\lambda})$ of the block $\R$ as follows:
\small
\begin{align*}
&V_{\R_\lambda} \doteq \{u_1, \ldots, u_m\} \cup \{\tilde u_1, \ldots, \tilde u_m\}
\cup \hspace{-5 mm} \bigcup\limits_{
\tiny
\begin{array}{c}
    i \in \boldsymbol{m}, \\
    d \in \boldsymbol{D_{\R_i}}
\end{array}}
\scriptsize
\hspace{-5 mm} \{x_{i,d}\}
\cup \hspace{-3 mm} \bigcup\limits_{
\tiny
\begin{array}{c}
    v \in \bar V, \\
    i \in \boldsymbol{m}
\end{array}}
\scriptsize
\hspace{-3 mm} \{f_{v,i}\},\\
&\forall i \in \boldsymbol{m}, \forall d \in \boldsymbol{D_{\R_i}}, (u_i, x_{i,d}) \in E_{\R_\lambda} \Leftrightarrow \gamma_{\R_i}(d) \neq 0,\\
&\forall i \in \boldsymbol{m}, \forall d \in \boldsymbol{D_{\R_i}}, \forall v \in \bar V, (f_{v,i}, x_{i,d}) \in E_{\R_\lambda} \Leftrightarrow \gamma_{v,i}(d) \neq 0,\\
&\forall i \in \boldsymbol{m}, \forall d_1, d_2 \in \boldsymbol{D_{\R_i}}, (x_{i,d_1}, x_{i,d_2}) \in E_{\R_\lambda} \Leftrightarrow d_1 = d_2 +1,\\
&\forall i \in \boldsymbol{m}, (x_{i,1}, \tilde u_i) \in E_{\R_\lambda}.
\end{align*}
\begin{proof}
By applying the principle of superposition to \eqref{eqPropGz} and \eqref{eqPropFz} we derive, for each component $i \in \boldsymbol{m}$ of the input signal, a state space representation
\begin{align*}
\dot{x_i}(t)&=A_i x(t)+B_i u_i(t)+ \sum\limits_{v \in \bar V} F_i f_{v,i}(t),\\
\tilde u_i(t) &= C_i x_i(t)
\end{align*}
of $\widetilde{U}_i(z) = \R_i(z) U_i(z) + \sum\limits_{v \in \bar V} \mathcal{F}_{v,i}(z) F_{v,i}(z)$. In particular, we define $x = [x_{i,1} \ldots x_{i,D_{\R_i}}]^{\top}$ and:
\begin{align}
& A_{i}=\left[
          \begin{array}{lll}
            0 & \vline & \textbf{I}_{D_{\R_i} -1}\\
            \hline
            0 & \vline & \textbf{0}_{1 \times (D_{\R_i} -1)}
          \end{array}
        \right], \notag\\
& B_{i}=\left[
          \begin{array}{ccc}
            \gamma_{\R_i}(1) & \cdots & \gamma_{\R_i}(D_{\R_i}) \\
          \end{array}
        \right]^{\top}, \notag\\
& F_{i}=\left[
          \begin{array}{ccc}
            \gamma_{v,i}(1) & \cdots & \gamma_{v,i}(D_{v,i}) \\
          \end{array}
        \right]^{\top}, \notag\\
& C_{i}=\left[
          \begin{array}{cccc}
            1 & 0 & \cdots & 0 \\
          \end{array}
        \right]. \label{eqStateSpaceSubsystem}
\end{align}
The result follows by definition of graph representation of a structured system.
\end{proof}
\end{proposition}

\normalsize

\begin{figure}[ht]
\begin{center}
\includegraphics[width=0.5\textwidth]{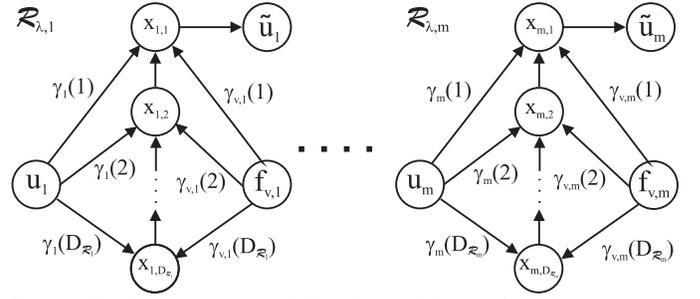}
\vspace{-0.8cm}
\caption{Graph representation of $\R_\lambda$ when a failure in the node $v$ occurs.}\label{figStructuredGraphWithFailures}
\vspace{-0.5cm}
\end{center}
\end{figure}

Figure \ref{figStructuredGraphWithFailures} provides an example of the graph representation of $\R_\lambda$ when a failure in the node $v$ occurs. The $\gamma$'s labeling some edges of $E_{\R_\lambda}$ just indicate that such edges are present if and only if the corresponding $\gamma$'s are not equal to 0. Note that $(V_{\R_\lambda}, E_{\R_\lambda})$ is composed by $m$ weakly connected components, each associated to the data flow of the $i$-th input component. Also note that each node $x_{i,d}$ is a variable associated to the $i$-th input component that will be delivered with a delay $d$ to the actuator node $v_{u,i}$. Finally, note that the sets of input and output nodes are respectively $U \doteq \{u_1, \ldots, u_m\}$ and $\tilde U \doteq \{\tilde u_1, \ldots, \tilde u_m\}$.

The same holds for defining the structured graph representation $(V_{\O_\lambda}, E_{\O_\lambda})$ of block $\O$, where the sets of input and output nodes are respectively $\tilde Y \doteq \{\tilde y_1, \ldots, \tilde y_\ell\}$ and $Y \doteq \{y_1, \ldots, y_\ell\}$.

Finally, let $(V_{\P_\lambda}, E_{\P_\lambda})$ be the structured graph representation of the plant $\P$, where the sets of input and output nodes are respectively $\tilde U \doteq \{\tilde u_1, \ldots, \tilde u_m\}$ and $\tilde Y \doteq \{\tilde y_1, \ldots, \tilde y_\ell\}$.

The model of a MCN $\M$ is the cascade of the blocks $\R$, $\P$ and $\O$. As a consequence, its structured graph representation $(V_{\M_\lambda}, E_{\M_\lambda})$ is given by the union of the structured graph representations of $\R_\lambda$, $\P_\lambda$ and $\O_\lambda$, and the set of nodes $\tilde U$ and $\tilde Y$ represent the interconnection nodes among such graphs. It is easy to see that all nodes in $\tilde U$ and $\tilde Y$ are (weak) bridges of $(V_{\M_\lambda}, E_{\M_\lambda})$ since the removal of one of them increases the number of weakly connected components.


\section{FDI of failures on structured systems} \label{secFDIstruct}

\begin{assumption}\label{assSimultaneousFailures}
We assume in this section that $f_{v,i}(k) = f_v(k)$ for all $i$ such that there exist $v' \in V_\R$, $h \in \boldsymbol{\Pi}$ with $(v,v') \in \eta_{\R_i}(h)$.
\end{assumption}
In other words, when a failure on a node $v$ occurs, it equally affects all the input components routed via $v$. This assumption is satisfied when $f_{v}(k)$ models a node failure, or when the malicious attack is not able to access separately the input components. We will discuss the case when this assumption does not hold in Section \ref{secRemoveAssumption1}.

Let a structured system $S_{\lambda}$ be given in the form:
\begin{align*}
\footnotesize
{x}(k+1) &= Ax(k) + Bu(k) + E_1 d(k) + F_1 f(k),\\
y(k) &= Cx(k) + Du(k) + E_2 d(k) + F_2 f(k),
\end{align*}
\normalsize
where $d(k)$ is a vector of disturbance signals and $f(k) = [f_1(k), \ldots, f_r(k)]^{\top}$ is a vector of $r$ failure signals. In this paper we will not consider the disturbance (i.e. $E_1=E_2=\boldsymbol{0}$) and leave such generalization to further work. In \cite{CommaultTAC2002} necessary and sufficient conditions have been derived on the graph representation of $S_\lambda$ that (generically) guarantee the existence of a bank of Luenberger observers, which takes as inputs $u(k), y(k)$, generates as output the \emph{residual} signals vector $\hat f(k) = [\hat f_1(k), \ldots, \hat f_r(k)]^{\top}$ and is characterized by a transfer function
\begin{equation*}
\footnotesize
\left[
  \begin{array}{c}
    \hat F_1(z) \\
    \vdots \\
    \hat F_r(z) \\
  \end{array}
\right] =
\left[
  \begin{array}{ccc}
    T_{11}(z) & \cdots &0 \\
    \vdots &  \ddots & \vdots \\
    0 & \cdots &  T_{rr}(z) \\
  \end{array}
\right]
\left[
  \begin{array}{c}
    F_1(z) \\
    \vdots \\
    F_r(z) \\
  \end{array}
\right],
\end{equation*}
\normalsize
where $\forall i \in \boldsymbol{r}$, $T_{ii}(z) \neq 0$. Characterizing the existence of such bank of observers is called the \emph{bank of observer-based diagonal FDI problem}. It is well known that the control input effects can be taken into account in the observer structure, therefore we will consider without loss of generality $B=D=0$. The theorem below characterizes the {bank of observer-based diagonal FDI problem} when there are no disturbances, and is a particular case of Theorem 3 in \cite{CommaultTAC2002}.
\begin{theorem}[from Theorem 3 in  \cite{CommaultTAC2002}]\label{thStructFDI}
The bank of observer-based diagonal FDI problem is generically solvable for a system $S$ if and only if ($i$) $S$ is structurally observable and ($ii$) $k = r$, where $k$ is the maximum number of fault-output vertex disjoint paths in the graph representation of $S_\lambda$.
\end{theorem}




\section{FDI of node failures on MCN\lowercase{s}} \label{secFDIMCN}

Given a MCN $\M$ subject to failures on communication nodes, if the {bank of observer-based diagonal FDI problem} is generically solvable for node failures in $\M_\lambda$, then the residual signals can be used to detect and identify possibly simultaneous occurrence of node failures.
In \cite{SmarraNecSys2012} we proved that, given a MIMO MCN $\M$ and if the plant $\P$ is controllable and observable, it is always possible to design a weight function $W$ such that $\M$ is controllable and observable. As a consequence, we assume in this paper that $\M_\lambda$ always satisfies Condition $(i)$ of Theorem \ref{thStructFDI}.

In order to design the network configuration of a MCN $\M$ to enable FDI of node failures, as the main result of this paper we state a formal relation between the network topology $G$ and scheduling/routing $\eta$ of $\M$, and solvability conditions of the bank of observer-based diagonal FDI problem for node failures. With this aim we first propose an algorithm to construct a graph $(\mathcal{V}, \mathcal{E})$, which essentially consists of the disjoint union of the controllability and observability graphs associated to each schedule, and of the structured graph representation of the plant. Then we will state a formal relation between solvability conditions of the {bank of observer-based diagonal FDI problem} for $\M_\lambda$ and $(\mathcal{V}, \mathcal{E})$.

\begin{figure}[ht]
\begin{center}
\includegraphics[width=0.5\textwidth]{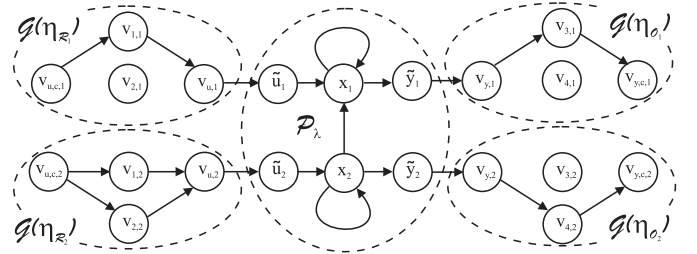}
\vspace{-0.6cm}
\caption{Graph $(\mathcal{V}, \mathcal{E})$ of Example \ref{exGraph4FDI}.}\label{figGraph4FDI}
\vspace{-0.6cm}
\end{center}
\end{figure}

\begin{definition}\label{algGraph}
Given a MCN $\M$, we define a graph
\scriptsize
\begin{equation}\label{eqAlgGraph1}
(\mathcal{V}, \mathcal{E}) = \left(\biguplus\limits_{i=1}^{m} G_\R(\eta_{\R_i})\right) \biguplus (V_{\P_\lambda}, E_{\P_\lambda}) \biguplus \left(\biguplus\limits_{i=1}^{\ell} G_\O(\eta_{\O_i})\right),
\end{equation}
\normalsize
with the addition of the following edges:
\begin{align}
\forall i \in \boldsymbol{m},\ &(v_{u,i}, \tilde u_i) \in \mathcal{E},\label{eqAlgGraph2}\\
\forall i \in \boldsymbol{\ell},\ &(\tilde y_i, v_{y,i}) \in \mathcal{E}.\label{eqAlgGraph3}
\end{align}
The disjoint union operator induces a map $\Gamma_\R$ that associates to each communication node $v \in V_\R$ a set of $m$ corresponding vertices $\{v_1, \ldots, v_m\} \subset \mathcal V$, one for each graph $G_\R(\eta_{\R_i}), i \in \boldsymbol m$. A map $\Gamma_\O$ can be defined similarly.
\end{definition}
\bigskip
\begin{example}\label{exGraph4FDI}
Consider a MCN $\M$ where the plant is given by $A = [1, 2; 0, 3]$, $B = C = \textbf{I}_{2}$, and $\scriptsize G_\R = \{v_{u,c}, v_1, v_2, v_{u,1}, v_{u,2}\}$, $G_\O = \{ v_{y,1}, v_{y,2}, v_3, v_4, v_{y,c}\},\ \eta_{\R_1} = \{(v_{u,c}, v_1), (v_1, v_{u,1})\}$, $\eta_{\O_1} = \{(v_{y,1}, v_3), (v_3, v_{y,c})\}, \eta_{\R_2} = \{(v_{u,c}, v_1), (v_{u,c}, v_2), (v_1, v_{u,2}), (v_2, v_{u,2})\}$, $\eta_{\O_2} = \{(v_{y,2}, v_4), (v_4, v_{y,c})\}$. \normalsize The corresponding graph $(\mathcal{V}, \mathcal{E})$ constructed as in Definition \ref{algGraph} is depicted in Figure \ref{figGraph4FDI}. It is easy to see that $\small \Gamma_\R(v_{u,c}) = \{v_{u,c,1}, v_{u,c,2}\}, \Gamma_\R(v_1) = \{v_{1,1}, v_{1,2}\}$ \normalsize etc. Note that the nodes $v_{u,1}, v_{u,2}, v_{y,1}, v_{y,2}$ do not split in the disjoint union, since each of them only belongs respectively to $G_\R(\eta_{\R_1})$, $G_\R(\eta_{\R_2})$, $G_\O(\eta_{\O_1})$, $G_\O(\eta_{\O_2})$.
\end{example}

\begin{lemma}\label{lemFDIonMCN}
Given a MCN $\M$ and a set $\{v_1, \ldots, v_r\} \subseteq V_\R \cup V_\O$ of faulty nodes, then the bank of observer-based diagonal FDI problem is generically solvable for node failures in $\M_\lambda$ if and only if there exists an $r$-linking from $\{f_{v_1}, \ldots, f_{v_r}\}$ to $Y$ in $(V_{\M_\lambda}, E_{\M_\lambda})$.
\begin{proof}
Straightforward by Theorem \ref{thStructFDI} and by definition of $\M_\lambda$.
\end{proof}
\end{lemma}

\begin{theorem}\label{thRelationFDIMCNConstructedGraph}
Let a MCN $\M$ and the associated graph $(\mathcal{V}, \mathcal{E})$ constructed as in Definition \ref{algGraph} be given. Given any set $\{v_1, \ldots, v_r\} \subseteq V_\R \cup V_\O$, there exists an $r$-linking from $\{f_{v_1}, \ldots, f_{v_r}\}$ to $Y$ in $(V_{\M_\lambda}, E_{\M_\lambda})$ if and only if there exists an $r$-linking from a set $\{\bar v_1 \in \Gamma(v_1), \ldots, \bar v_r \in \Gamma(v_r)\}$ to $\{v_{y,c,1}, \ldots, v_{y,c,\ell}\}$ in $(\mathcal{V}, \mathcal{E})$.

\begin{proof}
($\Rightarrow$) Consider a set of $r$ nodes given by the union of $\{v_{\R_1}, \ldots, v_{\R_p}\} \subseteq V_\R$ and $\{v_{\O_1}, \ldots, v_{\O_q}\} \subseteq V_\O$, with $p+q=r$, and define $F \doteq F_\R \cup F_\O$, where $F_\R = \{f_{v_{\R_1}}, \ldots, f_{v_{\R_p}}\}$ and $F_\O = \{f_{v_{\O_1}}, \ldots, f_{v_{\O_q}}\}$. If there exists an $r$-linking from $F$ to $Y$ in $(V_{\M_\lambda}, E_{\M_\lambda})$, then the following hold:
\begin{align}
&\exists \text{ a } p\text{-linking from } F_\R \text{ to a set } Y_\R = \{y_{{i_1}}, \ldots, y_{{i_p}}\} \subseteq Y, \label{eqProof1}\\
&\exists \text{ a } q\text{-linking from } F_\O \text{ to a set } Y_\O = \{y_{{j_1}}, \ldots, y_{{j_q}}\} \subseteq Y, \label{eqProof2}\\
&Y_\R \cap Y_\O = \varnothing, \label{eqProof3}
\end{align}
where $\{i_1,\ldots,i_p\} \subseteq \boldsymbol{\ell}$ and $\{j_1,\ldots,j_q\} \subseteq \boldsymbol{\ell}$. Since in $(V_{\M_\lambda}, E_{\M_\lambda})$ $\tilde U$ and $\tilde Y$ are bridge nodes and do not belong to any cycle, \eqref{eqProof1} implies the following:
\begin{align}
&\exists \text{ a } p\text{-linking from } F_\R \text{ to a set }\tilde U_\R = \{\tilde u_{{k_1}}, \ldots, \tilde u_{{k_p}}\} \subseteq \tilde U, \label{eqProof4}\\
&\exists \text{ a } p\text{-linking from } \tilde U_\R \text{ to a set }\tilde Y_\R = \{\tilde y_{{i_1}}, \ldots, \tilde y_{{i_p}}\} \subseteq \tilde Y, \label{eqProof5}\\
&\exists \text{ a } p\text{-linking from } \tilde Y_\R \text{ to }Y_\R, \label{eqProof6}
\end{align}
where $\{k_1,\ldots,k_p\} \subseteq \boldsymbol{m}$. By \eqref{eqProof4}, \eqref{eqAlgGraph1} and for the symmetric weakly connected components structure of $\R_{\lambda_i}$ and $G_\R(\eta_{\R_i})$, $i \in \boldsymbol{m}$, it follows that in $(\mathcal{V}, \mathcal{E})$
\begin{equation}
\exists \text{ a } p\text{-linking from }\{v_{\R_1, k_1}, \ldots, v_{\R_p, k_p}\} \text{ to }\{v_{u, k_1}, \ldots, v_{u, k_p}\}, \label{eqProof7}
\end{equation}
where $\forall \vartheta \in \boldsymbol{p}$, $v_{\R_\vartheta, k_\vartheta} \in \Gamma_\R(v_{\R_\vartheta})$. By \eqref{eqAlgGraph2} it follows that
\begin{equation}
\exists \text{ a } p\text{-linking from }\{v_{u, k_1}, \ldots, v_{u, k_p}\} \text{ to }\tilde U_\R. \label{eqProof8}
\end{equation}
By \eqref{eqProof5} and since \eqref{eqAlgGraph1} the structured graph representation $(V_{\P_\lambda}, E_{\P_\lambda})$ of the plant belongs both to $(\mathcal{V}, \mathcal{E})$ and $(V_{\M_\lambda}, E_{\M_\lambda})$, it follows that
\begin{equation}
\exists \text{ a } p\text{-linking from }\tilde U_\R \text{ to }\tilde Y_\R. \label{eqProof9}
\end{equation}
By \eqref{eqAlgGraph3} it follows that
\begin{equation}
\exists \text{ a } p\text{-linking from }\tilde Y_\R \text{ to }\{v_{y, i_1}, \ldots, v_{y, i_p}\}. \label{eqProof10}
\end{equation}
By \eqref{eqProof6}, \eqref{eqAlgGraph1} and for the symmetric weakly connected components structure of $\O_{\lambda_i}$ and $G_\O(\eta_{\O_i})$, $i \in \boldsymbol{\ell}$, it follows that
\begin{equation}
\exists \text{ a } p\text{-linking from }\{v_{y, i_1}, \ldots, v_{y, i_p}\} \text{ to }\{v_{y, c, i_1}, \ldots, v_{y, c, i_p}\}. \label{eqProof11}
\end{equation}
By \eqref{eqProof7}, \eqref{eqProof8}, \eqref{eqProof9}, \eqref{eqProof10} and \eqref{eqProof11} it follows that
\begin{align}
\exists \text{ a } p\text{-linking from }&\{v_{\R_1, k_1}, \ldots, v_{\R_p, k_p}\}\notag\\
\text{ to }&\{v_{y, c, i_1}, \ldots, v_{y, c, i_p}\}. \label{eqProof12}
\end{align}
By \eqref{eqProof2}, \eqref{eqAlgGraph1} and for the symmetric weakly connected components structure of $\O_{\lambda_i}$ and $G_\O(\eta_{\O_i})$, $i \in \boldsymbol{\ell}$, it follows that
\begin{align}
\exists \text{ a } q\text{-linking from }&\{v_{\O_1, j_1}, \ldots, v_{\O_q, j_q}\}\notag\\
\text{ to }&\{v_{y, c, j_1}, \ldots, v_{y, c, j_q}\}, \label{eqProof13}
\end{align}
where $\forall \vartheta \in \boldsymbol{q}$, $v_{\O_\vartheta, i_\vartheta} \in \Gamma_\O(v_{\O_\vartheta})$. By \eqref{eqProof3}, \eqref{eqProof12} and \eqref{eqProof13} it follows that there exists an $p+q$-linking in $(\mathcal{V}, \mathcal{E})$ from the $p+q$ dimensional set $\{v_{\R_1, i_1}, \ldots, v_{\R_p, i_p}, v_{\O_1, j_1}, \ldots, v_{\O_q, j_q}\}$ to the $p+q$ dimensional set $\{v_{y, c, i_1}, \ldots, v_{y, c, i_p}, v_{y, c, j_1}, \ldots, v_{y, c, j_q}\}$. Since $p+q=r$, this completes the proof.

($\Leftarrow$) Consider a set of $r$ nodes given by the union of $\{v_{\R_1}, \ldots, v_{\R_p}\} \subseteq V_\R$ and $\{v_{\O_1}, \ldots, v_{\O_q}\} \subseteq V_\O$, with $p+q=r$, and let $\mathcal F \doteq \mathcal F_\R \cup \mathcal F_\O$, where $\mathcal F_\R = \{v_{\R_1, k_1}, \ldots, v_{\R_p, k_p}\} \subseteq \mathcal V$ and $F_\O = \{v_{\O_1, j_1}, \ldots, v_{\O_q, j_q}\} \subseteq \mathcal V$, where $\{k_1,\ldots,k_p\} \subseteq \boldsymbol{m}$ and $\{j_1,\ldots,j_p\} \subseteq \boldsymbol{\ell}$. If there exists an $r$-linking from $\mathcal F$ to $\{v_{y,c,1}, \ldots, v_{y,c,\ell}\}$ in $(\mathcal{V}, \mathcal{E})$, then the following hold:
\begin{align}
&\exists \text{ a } p\text{-linking from } \mathcal F_\R \text{ to a set } \{v_{y,c,i_1}, \ldots, v_{y,c,i_p}\}, \label{eqRevProof1}\\
&\exists \text{ a } q\text{-linking from } \mathcal F_\O \text{ to a set } \{v_{y,c,j_1}, \ldots, v_{y,c,j_p}\}, \label{eqRevProof2}\\
&\{v_{y,c,i_1}, \ldots, v_{y,c,i_p}\} \cap \{v_{y,c,j_1}, \ldots, v_{y,c,j_p}\} = \varnothing, \label{eqRevProof3}
\end{align}
where $\{i_1,\ldots,i_p\} \subseteq \boldsymbol{\ell}$. By \eqref{eqAlgGraph2}, \eqref{eqAlgGraph3}, and since in $(\mathcal{V}, \mathcal{E})$ $\tilde U$ and $\tilde Y$ are bridge nodes that do not belong to any cycle, \eqref{eqRevProof1} implies the following:
\begin{align}
&\exists \text{ a } p\text{-linking from } \mathcal F_\R \text{ to a set }\tilde U_\R = \{\tilde u_{\R_{k_1}}, \ldots, \tilde u_{\R_{k_p}}\} \subseteq \tilde U, \label{eqRevProof4}\\
&\exists \text{ a } p\text{-linking from } \tilde U_\R \text{ to a set }\tilde Y_\R = \{\tilde y_{\R_{i_1}}, \ldots, \tilde y_{\R_{i_p}}\} \subseteq \tilde Y, \label{eqRevProof5}\\
&\exists \text{ a } p\text{-linking from } \tilde Y_\R \text{ to }\{v_{y,c,i_1}, \ldots, v_{y,c,i_p}\}. \label{eqRevProof6}
\end{align}
By \eqref{eqRevProof4}, \eqref{eqAlgGraph1} and for the symmetric weakly connected components structure of $\R_{\lambda_i}$ and $G_\R(\eta_{\R_i})$, $i \in \boldsymbol{m}$, it follows that in $(V_{\M_\lambda}, E_{\M_\lambda})$
\begin{equation}
\exists \text{ a } p\text{-linking from }\{v_{\R_1}, \ldots, v_{\R_p}\} \text{ to }\tilde U_\R. \label{eqRevProof7}
\end{equation}
Reasoning as above for \eqref{eqRevProof7}, by \eqref{eqRevProof5} and \eqref{eqRevProof6} it follows that
\begin{align}
&\exists \text{ a } p\text{-linking from }\tilde U_\R \text{ to }\tilde Y_\R, \label{eqRevProof9}\\
&\exists \text{ a } p\text{-linking from }\tilde Y_\R \text{ to } Y_\R = \{y_{{i_1}}, \ldots, y_{{i_p}}\} \subseteq Y. \label{eqRevProof11}
\end{align}
By \eqref{eqRevProof7}, \eqref{eqRevProof9} and \eqref{eqRevProof11} it follows that
\begin{equation}
\exists \text{ a } p\text{-linking from }\{v_{\R_1}, \ldots, v_{\R_p}\} \text{ to }Y_\R. \label{eqRevProof12}
\end{equation}
By \eqref{eqRevProof2}, \eqref{eqAlgGraph1} and for the symmetric weakly connected components structure of $\O_{\lambda_i}$ and $G_\O(\eta_{\O_i})$, $i \in \boldsymbol{\ell}$, it follows that
\begin{equation}
\exists \text{ a } q\text{-linking from }\{v_{\O_1}, \ldots, v_{\O_q}\} \text{ to }Y_\O = \{y_{{j_1}}, \ldots, y_{{j_p}}\} \subseteq Y. \label{eqRevProof13}
\end{equation}
By \eqref{eqRevProof3}, \eqref{eqRevProof12} and \eqref{eqRevProof13} it follows that there exists an $p+q$-linking in $(V_{\M_\lambda}, E_{\M_\lambda})$ from the $p+q$ dimensional set $\{v_{\R_1}, \ldots, v_{\R_p}, v_{\O_1}, \ldots, v_{\O_q}\}$ to the $p+q$ dimensional set $Y_\R \cup Y_\O \subseteq Y$. Since $p+q=r$, this completes the proof.
\end{proof}
\end{theorem}
The following corollary provides necessary and sufficient conditions on network topology, scheduling and routing for FDI of node failures over a MIMO MCN.
\begin{corollary}\label{corMainResult}
Given a MCN $\M$, the associated graph $(\mathcal{V}, \mathcal{E})$ constructed as in Definition \ref{algGraph} and a set $\{v_1, \ldots, v_r\} \subseteq V_\R \cup V_\O$ of faulty nodes, then the bank of observer-based diagonal FDI problem is generically solvable for node failures in $\M_\lambda$ if and only if there exists an $r$-linking from a set $\{\bar v_1 \in \Gamma(v_1), \ldots, \bar v_r \in \Gamma(v_r)\}$ to $\{v_{y,c,1}, \ldots, v_{y,c,\ell}\}$ in $(\mathcal{V}, \mathcal{E})$.
\begin{proof}
Straightforward by Lemma \ref{lemFDIonMCN} and Theorem \ref{thRelationFDIMCNConstructedGraph}.
\end{proof}
\end{corollary}
\begin{example}\label{exApplyResults}
Consider the MCN as in Example \ref{exGraph4FDI}. It is easy to see that the conditions of Corollary \ref{corMainResult} are satisfied for any 2-dimensional set of communication nodes, except for $(v_2, v_4)$. As a consequence, the FDI problem cannot be solved for 2 simultaneous failures. However, by adding to any time slot of the scheduling function $\eta_{\R_1}$ the transmission of links $(v_{u,c}, v_2)$ and $(v_2, v_{u,1})$, the conditions are satisfied and we can detect and isolate up to 2 simultaneous failures. An alternative solution is adding to the scheduling function $\eta_{\O_1}$ the transmission of links $(v_{y,1}, v_4)$ and $(v_4, v_{y,c})$. It is interesting that, in order to guarantee FDI of faulty nodes, the link scheduling order is irrelevant.
\end{example}
Example \ref{exApplyResults} shows that our results, for small graphs, can be used to enable FDI by manual designing the graph topology, scheduling and routing over the graph $(\mathcal{V}, \mathcal{E})$. For more complex networks graph theory algorithms can be exploited to automate the network design process, by searching for disjoint paths for any tentative set of faulty nodes in $2^{|V_\R| \cup |V_\O|}$, which is characterized by an exponential complexity. In the following we provide a sufficient condition that is easier to be verified.
\begin{lemma}\label{lemGraphConnectivity}
\cite{West01} Let a graph $(\mathcal V, \mathcal E)$ have connectivity $r$, and let $\mathcal V_1$, $\mathcal V_2$ be subsets of $\mathcal V$ each of size at least $r$, then there exists an $r$-linking from $\mathcal V_1$ to $\mathcal V_2$ (and vice versa).
\end{lemma}
\begin{proposition}\label{propFDIbyGraphConnectivity}
Given a MCN $\M$ and a positive integer $r \geq \min\{m,\ell\}$, let:
\begin{align}
&\forall v \in V_\R, |\Gamma_\R(v)| \geq r,\label{eqSuffCond1}\\
&(V_{\P_\lambda}, E_{\P_\lambda}) \text{ has connectivity } \geq r,\label{eqSuffCond2}\\
&\forall v \in V_\O, |\Gamma_\O(v)| \geq r.\label{eqSuffCond3}
\end{align}
Then for any $r$-dimensional set of faulty nodes the {bank of observer-based diagonal FDI problem} is generically solvable for $\M_\lambda$.
\begin{proof}
Consider a set of $r$ nodes given by the union of $\{v_{\R_1}, \ldots, v_{\R_p}\} \subseteq V_\R$ and $\{v_{\O_1}, \ldots, v_{\O_q}\} \subseteq V_\O$, with $p+q=r$. Let $(\mathcal{V}, \mathcal{E})$ be the graph constructed as in Definition \ref{algGraph} from $\M$, and let $\mathcal F \doteq \mathcal F_\R \cup \mathcal F_\O$, where $\mathcal F_\R = \{v_{\R_1, i_1}, \ldots, v_{\R_p, i_p}\} \subseteq \mathcal V$ and $\mathcal F_\O = \{v_{\O_1, j_1}, \ldots, v_{\O_q, j_q}\} \subseteq \mathcal V$, where $\{i_1,\ldots,i_p\} \subseteq \boldsymbol{\ell}$ and $\{j_1,\ldots,j_p\} \subseteq \boldsymbol{\ell}$. By \eqref{eqSuffCond1} it follows that there exists a $p$-linking from $\mathcal F_\R$ to a set $\{v_{u, k_1}, \ldots, v_{u, k_p}\}$, where $\{k_1,\ldots,k_p\} \subseteq \boldsymbol{m}$. By \eqref{eqAlgGraph2} and since $r \leq m$ it follows that there exists a $p$-linking from $\{v_{u, k_1}, \ldots, v_{u, k_p}\}$ to $\{\tilde u_{k_1}, \ldots, \tilde u_{k_p}\}$. By \eqref{eqSuffCond2} it follows that there exists a $p$-linking from $\{\tilde u_{k_1}, \ldots, \tilde u_{k_p}\}$ to $\{\tilde y_{i_1}, \ldots, \tilde y_{i_p}\}$. By \eqref{eqAlgGraph3} and since $r \leq \ell$ it follows that there exists a $p$-linking from $\{\tilde y_{i_1}, \ldots, \tilde y_{i_p}\}$ to a set $\{v_{y, i_1}, \ldots, v_{y, i_p}\}$. By \eqref{eqSuffCond3} it follows that there exists an $r$-linking from $\{v_{y, i_1}, \ldots, v_{y, i_p}\} \cup \mathcal F_\O$ to a set $\{v_{y, c, i_1}, \ldots, v_{y, c, i_p}, v_{y, c, j_1}, \ldots, v_{y, c, j_q}\}$. For all the above properties it follows that there exists an $r$-linking from $\mathcal F$ to $\{v_{y, c, i_1}, \ldots, v_{y, c, i_p}, v_{y, c, j_1}, \ldots, v_{y, c, j_q}\}$. By Corollary \ref{corMainResult}, and since $p+q=r$, this completes the proof.
\end{proof}
\end{proposition}
Conditions \eqref{eqSuffCond1} and \eqref{eqSuffCond3} can be verified in linear time with respect to $|V_\R| + |V_\O|$, but they are conservative. Indeed, it is easy to see that, for the solutions proposed in Example \ref{exApplyResults} where the FDI problem is solvable for up to 2 simultaneous failures, they are not satisfied. Condition \eqref{eqSuffCond2} can be verified by searching for disjoint paths on a graph characterized by cardinality $n+m+\ell$, and it is easy to show that it is a necessary condition (i.e. by assuming that all $r$ node failures occur in the controllability network).


\section{Removing Assumption \ref{assSimultaneousFailures}} \label{secRemoveAssumption1}

If we do not impose Assumption \ref{assSimultaneousFailures} a failure on a node $v \in V_\R$ affects all the input components routed via $v$ with possibly different signals $\{f_{v,i}(k)\}_{i \in \phi(v)}$, where
$$
\phi(v) \doteq \{i \in \boldsymbol{m} : (\exists v' \in V_\R, \exists h \in \boldsymbol{\Pi} : (v,v') \in \eta_{\R_i}(h))\}
$$
represents the set of input components routed via $v$.
In this case, since we aim at isolating node failures, we are just interested in detecting whether at least one of the signals $\{f_{v,i}(k)\}_{i \in \phi(v)}$ is active. As a consequence the conditions for solvability of the bank of observer-based diagonal FDI problem can be defined as follows.
\begin{lemma}\label{lemFDIonMCNnoAssumption1}
Given a MCN $\M$ and a set $\bar V = \{v_1, \ldots, v_r\} \subseteq V_\R \cup V_\O$ of faulty nodes, define $\bar F = \bigcup\limits_{v \in \bar V}\bigcup\limits_{i \in \phi(v)}\{f_{v,i}(k)\}$ the set of all failure signals. The bank of observer-based diagonal FDI problem is generically solvable for node failures in $\M_\lambda$ if and only if, for any $\bar v \in \bar V$ and any $\bar i \in \phi(\bar v)$, there exists an $r$-linking from $\bar F \setminus \bigcup\limits_{i \in \phi(\bar v) \setminus \{\bar i\}}\{f_{\bar v,i}(k)\}$ to $Y$ in $(V_{\M_\lambda}, E_{\M_\lambda})$.
\begin{proof}
By applying the same reasoning as in the proof of Theorem 3 in \cite{CommaultTAC2002} it can be stated that, given $\bar v \in \bar V$ and $\bar i \in \phi(\bar v)$, an $r$-linking exists from $\bar F \setminus \bigcup\limits_{i \in \phi(\bar v) \setminus \{\bar i\}}\{f_{\bar v,i}(k)\}$ to $Y$ in $(V_{\M_\lambda}, E_{\M_\lambda})$ if and only if it is possible to design an observer that generates a signature $\hat f_{\bar v, \bar i}(k)$ characterized by a non-zero transfer function with respect to $f_{\bar v, \bar i}(k)$, a zero transfer function with respect to any failure signal in the set $\bar F \setminus \bigcup\limits_{i \in \phi(\bar v) \setminus \{\bar i\}}\{f_{\bar v,i}(k)\}$, and a (possibly zero or non-zero) transfer function with respect to any failure signal in the set $\bigcup\limits_{i \in \phi(\bar v) \setminus \{\bar i\}}\{f_{\bar v,i}(k)\}$. Therefore, the signature $\hat f_{\bar v, \bar i}(k)$ is not affected by the failure signals of any other node $v \neq \bar v$, and is (generically) non-zero if $f_{\bar v, \bar i}(k)$ is non-zero. As a consequence, we can detect and isolate a failure acting on node $\bar v$ when at least one signature $\hat f_{\bar v, i}(k), i \in \phi(\bar v)$ is non-zero. This completes the proof.
\end{proof}
\end{lemma}
\begin{proposition}\label{propFDIonMCNnoAssumption1}
Given a MCN $\M$ and 2 faulty nodes $v_1, v_2 \in V_\R$, the bank of observer-based diagonal FDI problem is generically solvable for node failures in $\M_\lambda$ only if $\phi(v_1) \cap \phi(v_2) = \varnothing$.
\begin{proof}
Assume that $\phi(v_1) \cap \phi(v_2) = \{k\} \in \boldsymbol{m}$, and let $\bar v = v_1$, $\bar i = k$. In the graph $(V_{\M_\lambda}, E_{\M_\lambda})$ there does not exist a 2-linking from $\{f_{v_1,k}, f_{v_2,k}\}$ to $Y$ since they only have outgoing links to the weakly connected component $\R_{\lambda, k}$ of the structured graph representation of the block $\R$ (see Figure \ref{figStructuredGraphWithFailures}). As a consequence, the maximal linking from $\bar F \setminus \bigcup\limits_{i \in \phi(v_1) \setminus \{k\}}\{f_{v_1,i}(k)\}$ to $Y$ in $(V_{\M_\lambda}, E_{\M_\lambda})$ is $k=1$. Since $r=2$ the results follows.
\end{proof}
\end{proposition}
The above proposition can be proven similarly for the observability graph and shows that, if Assumption \ref{assSimultaneousFailures} is not raised, the bank of observer-based diagonal FDI problem is generically solvable for node failures in $\M_\lambda$ only for trivial network topologies, namely when all graphs $G_\R(\eta_{\R_i}), i \in \boldsymbol{m}$ and $G_\O(\eta_{\O_i}), i \in \boldsymbol{\ell}$ consist of a single communication node (namely, they are not multi-hop networks). This can be easily seen by considering that, for any 2 nodes $v_1, v_2$ belonging to a given $G_\R(\eta_{\R_k})$, $k \in \phi(v_1) \cap \phi(v_2)$. This negative result is intuitive: if a malicious attack is able to access separately the input components and to inject unrelated signals to each of them, it is much more difficult to exploit redundancy to perform FDI. To overcome this difficulty an interesting venue for future work is providing milder conditions that enable to detect and isolate failures of clusters of nodes, instead of isolating singleton node failures. In particular, it would be interesting to compute for a given MCN the minimal node clustering that enables FDI.


%

\bibliographystyle{plain}
\bibliography{mcnbib}

\end{document}